\documentclass[11pt]{article}
\usepackage[abbr,agsmcite]{harvard}
\usepackage{algorithm}
\usepackage{algorithmic}
\usepackage[T1]{fontenc}              
\usepackage[latin9]{inputenc}         

\usepackage[a4paper]{geometry}        
\usepackage{amsmath,amssymb,amsthm}          
\usepackage{bm}                       

\usepackage{url}                      
\usepackage{xspace}                   

\usepackage{graphicx}

\usepackage{verbatim}  

\usepackage{setspace}                 
\newtheorem{algo}{Algorithm}
\newtheorem{corollary}[algo]{Corollary}

\newtheorem{lemma}[algo]{Lemma}
\newtheorem{proposition}[algo]{Proposition}
\newtheorem{thm}[algo]{Theorem}

\newcommand{\PP}{\mathcal{P}}

\bmdefine\g{g} \bmdefine\K{K}\bmdefine\X{X} \bmdefine\X{X}
\bmdefine\D{D} \bmdefine\K{K} \bmdefine\Z{Z} \bmdefine\x{x}
\bmdefine\z{z} \bmdefine\y{y} \bmdefine\Y{Y}
\bmdefine\bfalpha{\alpha} \bmdefine\bfmu{\mu} \bmdefine\M{M}
\bmdefine\Q{Q} \bmdefine\P{P} \bmdefine\w{w} \bmdefine\W{W}
\bmdefine\p{p} \bmdefine\T{T} \bmdefine\t{t} \bmdefine\r{r}
\bmdefine\B{B} \bmdefine\I{I} \bmdefine\u{u} \bmdefine\p{p}
\bmdefine\Sig{\Sigma} \bmdefine\E{E} \bmdefine\F{F} \bmdefine\S{S}
\bmdefine\s{s} \bmdefine\w{w} \bmdefine\b{b} \bmdefine\W{W}
\bmdefine\w{w} \bmdefine\V{V} \bmdefine\v{v} \bmdefine\q{q}
\bmdefine\R{R} \bmdefine\A{A} \bmdefine\bflambda{\lambda}
\bmdefine\C{C} \bmdefine\U{U}
\bmdefine\d{d}

\newcommand{\bbeta}{\mbox{\boldmath$\beta$}}
\newcommand{\bLambda}{\mbox{\boldmath$\Lambda$}}
\newcommand{\bvarepsilon}{\mbox{\boldmath$\epsilon$}}

\newcommand{\bDelta}{\mbox{\boldmath$\Delta$}}

\newcommand{\balpha}{\mbox{\boldmath$\alpha$}}
\newcommand{\bTheta}{\mbox{\boldmath$\Theta$}}
\newcommand{\blambda}{\mbox{\boldmath$\lambda$}}

\title{Penalized Partial Least Squares Based on B-Splines Transformations}

\author{Nicole Kr\"amer \thanks{Department of Electrical Engineering and Computer
Science,
       Technical University of Berlin,
       Franklinstr. 28/29, 10587 Berlin,
Germany,\texttt{nkraemer@cs.tu-berlin.de}} \and Anne-Laure
Boulesteix\thanks{Department of Medical Statistics and Epidemiology,
Technical University of Munich, Ismaningerstr. 22, 81675 Munich, Germany,
\texttt{anne-laure.boulesteix@tum.de}} \and  Gerhard Tutz\thanks{Department
of Statistics, University of Munich, Akademiestr. 1, 80799 Munich, Germany, \texttt{tutz@stat.uni-muenchen.de}}}

\begin{document}
\maketitle

\begin{abstract}
We propose a novel method to model nonlinear regression problems by adapting
the principle of penalization to Partial Least Squares (PLS). Starting  with a
generalized additive model, we expand the additive component of each
variable in terms of a generous amount of B-Splines basis functions. In
order to prevent overfitting and to obtain smooth functions, we
estimate the regression model by applying a penalized version of
PLS. Although our motivation for penalized PLS stems from its use for
B-Splines transformed data, the proposed approach is very general and can be
applied to other penalty terms or to other dimension reduction techniques. It
turns out that penalized PLS can be computed virtually as fast
as PLS. We prove a close connection of penalized PLS to the solutions of
preconditioned linear systems. In the case of high-dimensional data, the new method is shown to be an attractive
competitor to other techniques for estimating generalized additive
models. If the number of predictor variables is high compared to the number
of examples, traditional  techniques often suffer from overfitting. We
illustrate that penalized PLS performs well in these situations.

\end{abstract}
\newpage
\section{Introduction}
Nonlinear regression effects may be modeled  via additive regression models of the form 
\begin{eqnarray}
\label{eq:model}
Y&=&\beta_0+f_1(X_1)+\dots+f_p(X_p) +\varepsilon\,.
\end{eqnarray}
where the functions $f_1,\dots,f_p$ have unspecified functional form.
An approach which allows flexible representation of the functions $f_1,\dots,f_p$
is the expansion in basis functions \cite{Hastie9001}. To prevent
overfitting, there are two general approaches. In the first approach, each function $f_j$ is the sum of only a small set of basis
functions,
\begin{eqnarray}
\label{eq:expansion}
f_j(x)&=&\sum_{k=1} ^{K_j} \beta_{kj} B_{kj}(x)\,. 
\end{eqnarray}
The basis functions $B_{kj}$ are  chosen adaptively by a selection procedure. The second
approach (that is outlined in Section \ref{sec:splines}) circumvents the
problem of basis function selection. Instead, we allow a generous amount
$K_j \gg1$ of
basis  functions in the expansion (\ref{eq:expansion}). As this usually leads
to high-dimensional and highly correlated data, we penalize  the coefficients
$\beta_{jk}$ in the estimation process \cite{Eilers9601}.

Quite generally, a different approach to deal with high dimensionality is to
use dimension reduction techniques such as Partial Least Squares (PLS) \cite{Wold7501,Wold8401}. The main
idea is to build a few components from the predictor variables and to regress $\y$ onto
these components.  A short overview on PLS can be found in Section \ref{sec:pls}.

As a linear approach, PLS probably fails to yield high prediction accuracy in
the case of nonlinear relationships between predictors and responses as in
(\ref{eq:model}). In order to incorporate nonlinear structures, it might be
advisable to 
transform the original predictors  preliminarily to a PLS regression.  This approach has been proposed by \citeasnoun{DS1997}
and \citeasnoun{Dur2001} in different variants.
The method proposed by \citeasnoun{DS1997} is based on a variant of PLS that may
be computed via an iterative algorithm. They suggest an approach that
incorporates splines transformations of the predictors within each iteration of
the iterative algorithm. In contrast, the method proposed by \citeasnoun{Dur2001} is
global. The predictors are first transformed using splines basis functions as a
preliminary step, then PLS regression is performed on the transformed data
matrix. The choice of the degree $d$ of the polynomial pieces and of the number
of knots is performed  by an either ascending or descending search
procedure that is not automatic.

For large numbers of variables, this search procedure is
computationally intensive and might overfit the training data.
In the present article, we suggest an alternative approach based on the
penalty strategy of \citeasnoun{Eilers9601}. As described in Section \ref{sec:splines}, we transform the initial data
matrix  nonlinearly using  B-splines basis functions. Our new
method, which we call penalized PLS, is based on the following principle.
The equivalent of penalizing the (higher order) differences of adjacent
B-splines coefficients is, in the framework of dimension
reduction, the  penalization of (higher order) differences of adjacent
weights. 

In Section \ref{sec:penPLS}, we introduce an adaptation of the principle of
penalization  to PLS. More precisely, we present a penalized version of the
optimization problem attached to PLS. Although the motivation stems from its use for
B-splines transformed data, the proposed approach is very general and can be
adapted to other penalty terms or to other dimension reduction techniques
such as Principal Components Analysis. It turns out that the new method
shares a lot of properties of PLS and that its computation requires
virtually no  extra costs. We highlighten the close connection between penalized
PLS and preconditioned linear systems. It is already known that PLS is
equivalent to the conjugate gradient method \cite{Hestenes5201} applied to the set of normal
equations associated to a linear regression problem. We prove that penalized PLS corresponds to a conjugate gradient
method for  a preconditioned set of normal equations, where the
preconditioner depends on the penalty term. Furthermore, we show that this new
technique is closely related to the so-called kernel trick. More precisely, we prove 
that penalized PLS is equivalent to ordinary PLS using a generalized inner
product that is defined by the penalty term. In Sections \ref{sec:birth} and
\ref{sec:polymer}, we illustrate our
method on different data sets.

In the rest of the paper, we restrict ourselves to a univariate response. In Section \ref{sec:conclusion}, we stress that the extension of
our method to a multivariate response is straightforward. 
\section{Partial Least Squares Regression}
\label{sec:pls}
Let us consider the general linear regression problem. We want to
predict a univariate response variable $Y$ using $p$ predictor variables
$X_1,\ldots,X_p$ based on a finite set
\begin{displaymath}
\left\{\left(y_i,\x_i\right)=\left(y_i,x_{i1},\ldots,x_{ip}\right)
\,,\,i=1,\ldots,n\right \}
\end{displaymath}
of observations. We set
\begin{displaymath}
\X  =  \left(\begin{array}{c} \x_1^T \\ \dots \\ \x_n^T
\end{array}\right) \in \mathbb{R}^{n\times p}, \ \
\y  = \left(\begin{array}{c} y_1 \\ \dots \\ y_n \end{array}\right) \in \mathbb{R}^n\,,
\end{displaymath}
and require for simplicity of notation that both $\X$ and $\y$ are centered.
If we assume that
the relationship between predictors and response is linear, this
relationship can be represented in compact form by
\begin{eqnarray*}
\label{eq:linreg}
\y&=&\X\bbeta +\bvarepsilon\,\,.
\end{eqnarray*}
Here, $\bbeta$ is the $p$-dimensional vector of regression coefficients and
$\bvarepsilon$ is the vector of residuals.\\

When $n<p$, the usual regression tools such as ordinary least squares (OLS)
regression  cannot be applied to estimate $\bbeta$ since the $p\times p$ covariance matrix $(1/n)\X^T\X$ (which
has rank at most $n-1$) is singular.   From a technical point of view, this
may be solved by replacing the inverse of the covariance matrix by a generalized
inverse. However, for $n<p$, OLS usually fits the training data perfectly
and one cannot expect
the method to perform well on a new data set. Partial Least Squares (PLS)
\cite{Wold7501,Wold8401} is an alternative regression tool
which is more appropriate in the case of highly correlated predictors
and high-dimensional data. PLS is a standard tool for analyzing chemical data
\cite{Martens8901}, and in recent years, the success of PLS has lead to
applications in other scientific fields such as physiology
\cite{Rosipal0301} or bioinformatics \cite{Boulesteix0601}.

The main idea of PLS is to build orthogonal components $\t_1,\ldots,
\t_m$ from the original predictors $\X$  and to  use them as
predictors in a least squares regression. There are different PLS
techniques to extract these components, and each of them gives rise
to a different variant of PLS. It is not our aim to explain all
variants and we  focus on two of them. An overview on different
forms of PLS can be found in \citeasnoun{Rosipal0601}.
 A component is a linear combination
of the original predictors that hopefully reflects the relevant
structure of the data. PLS is similar to Principal Components
Regression (PCR). The difference is that PCR extracts components
that explain the variance in the predictor variables  whereas PLS
extracts components that have a large covariance with $\y$. We now
formalize this concept. A latent component $\t$ is a linear
combination $\t=\X \w$  of the predictor variables. The vector $\w$
is usually called the weight vector. We want to find a component
with maximal covariance to $\y$, that is we want to maximize the
empirical squared covariance
\begin{eqnarray*}
\text{cov}^2\left(\X \w,\y\right) &=&  \w ^T\X^T\y\y^T\X\w\,.
\end{eqnarray*}
We have to constrain $\w$ in order to obtain identifiability,
choosing
\begin{eqnarray}
\label{eq:crit1}\max && \w ^T \X^T \y \y^T \X \w\,,\\
\label{eq:constr1}\text{subject to}&& \|\w\|=1\,.
\end{eqnarray}
Using Lagrangian multipliers, we conclude that the solution $\w_1$ is -- up
to a scaling factor -- equal to  $\X^T \y$.

Let us remark that (\ref{eq:crit1}) and (\ref{eq:constr1}) are equivalent to
\begin{eqnarray}
\label{eq:critneu}
\max && \frac{\w^T \X^T \y \y^T \X \w}{\w^T \w}\,\,.
\end{eqnarray}
The solution of (\ref{eq:critneu}) is only unique up to a scalar. The
normalization of the weight vectors $\w$ to length $1$ is not essential
for the PLS algorithm and PLS algorithms differ in the way they scale
the weight vectors and components. In this paper, we present all algorithms
without the scaling of the vectors, in order to keep the notation as simple as possible.

Subsequent components $\t_2,\t_3,\ldots$ are chosen such that they
maximize (\ref{eq:crit1}) and that all components $\t_i$ are
mutually orthogonal. In PLS, there are different techniques to
extract subsequent components, and each technique gives rise to a
variant of PLS. We briefly introduce two of them. In the method
called SIMPLS \cite{deJong9301}, one computes  for the $i$th
component,
\begin{eqnarray*}
\max && \w ^T \X^T \y \y^T \X \w\,,\\
\text{subject to}&& \|\w\|=1 \text{ and } \X \w \perp \t_j ,j<i\,.
\end{eqnarray*}
Alternatively, one can deflate the original predictor variables
$\X$. That is, we only consider the part of $\X$ that is orthogonal
onto all components $\t_j,j<i$. For any matrix $\V$, let us denote
by $\mathcal{P}_{\V}$ the orthogonal projection onto the space that
is spanned by the columns of $\V$. In matrix notation, we have
\begin{eqnarray}
\label{eq:proj}
\mathcal{P}_{\V}&=& \V \left(\V^T \V\right)^{-} \V^T\,.
\end{eqnarray}
The deflation of $\X$ with respect to the components $\t_1,\ldots,\t_{i-1}$
is defined as
\begin{eqnarray}
\label{eq:deflation}
\X_i = \X - \mathcal{P}_{\t_1,\ldots,\t_{i-1}} \X = \X_{i-1} -
\PP_{\t_{i-1}} \X_{i-1}\,.
\end{eqnarray}
For the computation of the $i$th component, $\X$ is replaced by
$\X_i$ in (\ref{eq:crit1}). This method is called the NIPALS
algorithm \cite{Wold7501}. The two methods are equivalent if $\y$ is
univariate in the sense that we end up with the same components
$\t_i\,$ \cite{deJong9301}. In this paper, we use the NIPALS
algorithm. In summary, the  PLS algorithm is described in algorithm
\ref{algo:pls1}.
\begin{algorithm}
\caption{NIPALS algorithm}
\begin{algorithmic}
\label{algo:pls1} \STATE{Input: $X_1=X$, $\y$, number of components
$m$} \FOR{i=1,\ldots,m} \STATE{$\w_i= \X_i ^T \y$ (weight vector)}
\STATE{ $\t_i= \X_i \w_i$(component)} \STATE{ $\X_{i+1}=\X_i -
\mathcal{P}_{\t_i} \X_i$(deflation)} \ENDFOR
\end{algorithmic}
\end{algorithm}

PLS used to be overlooked by statisticians and was considered an
algorithm rather than a sound statistical model. This attitude is in parts 
understandable, as  in the early literature on the subject, PLS was
explained solely in terms of formulas as in algorithm \ref{algo:pls1}.  Due to its success in
applications, the interest in the statistical properties of PLS has
risen. It can be related to other dimension reduction techniques such as
Principal Components Regression and Ridge Regression and these methods can
be cast under a unifying framework \cite{Stone9001}. The shrinkage
properties of PLS have been studied extensively
\cite{Frank9301,deJong9501,Goutis9601,Butler0001}. Furthermore, it can be shown that PLS is
closely connected to Krylov subspaces and the conjugate gradient method \cite{Helland8801,Phatak0301}. We
discuss this method in more detail in Section \ref{sec:penPLS}.

Let us return to the PLS algorithm. With
\begin{eqnarray*}
\T&=& \left(\t_1,\ldots,\t_m\right)\,.
\end{eqnarray*}
denoting the collection of components, the fitted response is given
by
\begin{eqnarray}
\label{eq:Yhat} \widehat{\y}&=&\T(\T^T\T)^{-1}\T^T\y =
\mathcal{P}_{\T} \y \,.
\end{eqnarray}
In order to obtain the response for new observations, we have to
determine the vector of regression coefficients $\widehat \y=  \X
\widehat \bbeta\,$. Therefore, a representation of the components
$\t_i=\X_i \w_i$ as a linear combination of the original predictors
$\X$ is needed. In other words, we have to derive weight vectors
$\widetilde \w_i$ with
\begin{displaymath}
\X \widetilde \w_i = \X_i \w_i\,.
\end{displaymath}
They are in general different from the ``pseudo'' weight vectors $\w_i$ that
are computed by the NIPALS algorithm. In order to
avoid redundancy, the derivation of these weight vectors is deferred until
Section \ref{sec:penPLS}.

It should be noted that the number $m$ of PLS components is an
additional model parameter that has to be estimated. One way of
determining $m$ is by cross-validation.
\section{Penalized Regression Splines}
\label{sec:splines}
The fitting of generalized additive models by use of  penalized
regression splines has become a widely used tool in statistics.
Starting with the seminal paper by \citeasnoun{Eilers9601}, the approach
has been extended and applied in various publications 
\cite{Ruppert0201,Wood0001,Wood0601}. The basic concept is to expand the
additive component of each variable $X_j$ in basis functions as in
(\ref{eq:expansion}) and to estimate the coefficients by
penalization techniques. As suggested in \citeasnoun{Eilers9601},
B-splines  are used as basis functions yielding so-called
P-splines (for penalized B-splines). 
Splines are one-dimensional piecewise polynomial functions. The
points at which the pieces are connected are called knots or
breakpoints. We say that a spline is of order $d$ if all polynomials are of
degree $\leq d$ and if the spline is $(d-1)$ times continuously
differentiable at the breakpoints. A particular efficient set of basis
functions are B-splines \cite{deBoor7801}. The number of basis functions depends on the order of the splines and the
number of breakpoints. For  a given variable $X_j$, we consider a set of  corresponding B-splines basis functions
$B_{1j},\dots,B_{Kj}$. These basis functions define a nonlinear map
\begin{eqnarray*}
\Phi_j(x)&= &\left(B_{1j}(x),\ldots,B_{Kj}(x)\right)^T\,.
\end{eqnarray*}
By performing such a transformation on each of the variables
$X_1,\dots,X_p$, the observation vector $\x_i$  turns into a
vector
\begin{eqnarray}
\label{eq:Z}
\z_i&=&(B_{11}(x_{i1}),\dots,B_{m1}(x_{i1}),\dots,B_{1p}(x_{ip}),\dots,B_{mp}(x_{ip}))^T\\
\nonumber &=& \Phi(\x_i)
\end{eqnarray}
of length $pK$. Here $\Phi$  is the function defined by the B-splines. The resulting data matrix obtained
by the transformation of $\X$ has
dimensions $n\times pK$ and will be denoted by $\Z$ in the rest of
the paper. In the examples in Sections \ref{sec:birth} and
\ref{sec:polymer},  we consider the most widely used cubic
B-splines, i.e. we choose $d=3$.

The estimation of (\ref{eq:model}) is transformed into the estimation of the
$pK$-dimensional vector that consists of the coefficients $\beta_{jk}$:
\begin{displaymath}
\bbeta^T=\left(\beta_{11},\ldots,\beta_{K1},\ldots
\beta_{12},\ldots,\beta_{Kp}\right)= \left(\bbeta_{(1)} ^T,\ldots, \bbeta_{(p)} ^T\right)\,.
\end{displaymath}
As explained above, the vector $\bbeta$ determines a nonlinear, additive
function
\begin{displaymath}
f(\x)=\beta_0 + \sum_{j=1} ^p f_j(x_j)= \beta_0 +\sum_{j=1} ^p
\sum_{k=1} ^K \beta_{kj} B_{kj}(x_j)= \beta_0 + \Phi(\x)^T \bbeta\,.
\end{displaymath}
As $\Z$ is usually high-dimensional, the estimation of $\bbeta$ by
minimizing the squared error 
\begin{displaymath}
\frac{1}{n}\sum_{i=1} ^n \left( y_i - f(\x_i)\right)^2=
\frac{1}{n}\|\y -\beta_0 - \Z \bbeta\|^2
\end{displaymath}
usually leads to overfitting.
Following \citeasnoun{Eilers9601}, we use for each variable many basis
functions, say $K\approx 20$, and estimate by penalization. The idea is to
penalize the second derivative of the function $f$. \citeasnoun{Eilers9601} show
that the following difference penalty term is a good approximation of the
penalty on the second derivative of $f$,
\begin{eqnarray*}
P(\bbeta)&=&\sum_{j=1}^p\sum_{k=3}^m\lambda_j(\Delta^2
\beta_{kj})^2\,.
\end{eqnarray*}
These are also called the second-order differences of adjacent parameters. The difference operator $\Delta^2\beta_{kj}$ has the
form
\begin{eqnarray*}
\Delta^2\beta_{kj} & = & (\beta_{kj}-\beta_{k-1,j})-(\beta_{k-1,j}-\beta_{k-2,j})\\
& = & \beta_{kj}-2\beta_{k-1,j}+\beta_{k-2,j}.
\end{eqnarray*}
The coefficients $\lambda_j\geq0$ control the amount of penalization. This penalty term  can be expressed in terms of a
penalty matrix $\P$.
We denote by  $\D_K$ the $(K-1)\times K$ matrix
\[
\D_K\ \ =\ \
\left(
\begin{array}{rrrrr}
1 & -1 & . & .&  .\\
.& 1 & -1 & . &  .\\
.& .&. & .& . \\
.& .& .& 1 & -1
\end{array}
\right)
\]
that defines the first order difference operator. Setting
\begin{displaymath}
\K_2 =  (\D_{K-1}\D_K)^T\D_{K-1}\D_K\,,
\end{displaymath}
we conclude that the penalty term equals
\begin{displaymath}
P(\bbeta)= \sum_{j=1}^p \lambda_j \bbeta_{(j)} ^T \K_2
\bbeta_{(j)}= \bbeta^T (\bDelta_{\blambda}\otimes\K_2) \bbeta\,.
\end{displaymath}
Here  $\bDelta_{\blambda}$ is  the $p\times p$ diagonal matrix containing $\lambda_1,\dots,\lambda_p$
on its diagonal and $\otimes$ is the Kronecker product. The generalization of this method to higher-order differences of the
coefficients of adjacent B-splines is
straightforward.  We simply replace $\K_2$ by
\begin{displaymath}
\K_q=(\D_{K-q+1}\dots\D_K)^T(\D_{K-q+1}\dots\D_K)\,.
\end{displaymath}
To summarize, the penalized least squares criterion has the form
\begin{eqnarray}
\label{eq:gam}
\widehat R_{P} (\bbeta)&=&\frac{1}{n}\| \y - \beta_0 - \Z \bbeta\|^2 + \bbeta^T \P \bbeta
\end{eqnarray}
with the penalty matrix $\P$ defined as
\begin{eqnarray}
\label{eq:P}
\P&=& \bDelta_{\blambda} \otimes \K_q\,.
\end{eqnarray}
This is a symmetric matrix that is positive semidefinite.
\section{Penalized Partial Least Squares Regression}
\label{sec:penPLS}
We now introduce a general framework to combine PLS with penalization
terms. We remark that this is not limited to spline transformed variables or
to the special shape of the penalty matrix $\P$ that is defined in
(\ref{eq:P}). For this reason, we present the new method in terms of the
original data matrix $\X$ and only demand that $\P$  is a symmetric matrix such that $\I_p+\P$ is positive definite.

Again, we restrict ourselves to univariate responses $\y$.  Penalized PLS for
multivariate responses is briefly discussed in Section
\ref{sec:conclusion}. We modify the optimization criterion (\ref{eq:critneu}) of PLS in the
following way.  The first component $\t_1=\X\w_1$ is  defined by the
solution  of the problem
\begin{eqnarray}
\label{eq:critpen}
\text{arg}\max_{\w} \frac{\w^T \X^T \y \y^T \X\w}{\w^T \w +  \w^T \P \w}\,.
\end{eqnarray}
Using Lagrangian multipliers, we obtain  the solution
\begin{eqnarray}
\label{eq:penweight} \w_1&=& \M\X^T \y
\end{eqnarray}
with $\M= \left(\I_p +\P\right)^{-1}\,$. Subsequent weight vectors
and components are computed by deflating $\X$ as described in
(\ref{eq:deflation})  and then maximizing (\ref{eq:critpen}) with
$\X$ replaced by $\X_i$. In particular, we can compute the weight
vectors and components of penalized PLS by simply replacing $\w_i=
\X_i ^T \y$ by (\ref{eq:penweight}) in algorithm \ref{algo:pls1}.

We now present  results on penalized PLS that allow us to compute its regression
vectors efficiently. Note that all results on penalized PLS also hold for
ordinary PLS if we choose  $\P={\bf 0}$. Let
\[
\begin{array}{rclcrcl}
\T&=&\left(\t_1,\ldots,\t_m\right)&,&\W&=& \left(\w_1,\ldots,\w_m\right) \,,
\end{array}
\]
denote  the matrices of components and weight vectors respectively.
\begin{lemma}\label{lemma:R}
The matrix
\begin{eqnarray*}
\R&=& \T^T \X\ \W\in \mathbb{R}^{m\times m}
\end{eqnarray*}
is upper bidiagonal, that is
\begin{displaymath}
r_{ij}= \t_i ^T \X \w_j = 0
\end{displaymath}
if $i<j$ or $i+1>j$. The matrix $\R$ is invertible. Furthermore, the columns
of $\T$ and the columns of $\X\W$ span the same space.
\end{lemma}
This is an extension of a result for ordinary PLS that can be found
e.g. in \citeasnoun{Manne8701}. The proof can be found in the appendix. We
can now determine the regression coefficients for penalized PLS.
\begin{proposition}
\label{regvector}
The Penalized PLS regression vector obtained after $m$ steps is
\begin{eqnarray}
\label{eq:beta}
\widehat\bbeta^{(m)} _{PPLS}&=& \W\left(\W^T \X^T \X \W\right)^{-1} \W^T \X^T \y\,.
\end{eqnarray}
In particular, the penalized PLS estimator is the solution of the
constrained minimization problem
\begin{eqnarray}
 \min_{\bbeta}& &\|\y-\X\bbeta\|^2  \nonumber\\
\label{eq:restr.ols}\text{subject to}& &\bbeta \in \text{span}\left\{\w_1,\ldots,\w_m\right\}.
\end{eqnarray}
\end{proposition}
\begin{proof}
 We deduce from lemma \ref{lemma:R} that the columns of $\X \w$ span the same
space as the columns of $\T$. As  PLS is ordinary least squares regression with predictors $\t_1,\ldots,\t_m$, we have
\begin{displaymath}
\widehat \y= \PP_{\T} \y= \PP_{\X\W} \y= \X \W \left( \W^T \X^T \X \W
\right)^{-1} \W^T \X^T \y\,.
\end{displaymath}
The second statement can be proven by noting that the OLS
minimization problem with constraints (\ref{eq:restr.ols}) is equivalent to an unconstrained minimization problem
for $\bbeta= \W \balpha$ with $\balpha \in \mathbb{R}^m$. If we
plug this into the formula for the OLS estimator, we obtain (\ref{eq:beta}).
\end{proof}
Formula (\ref{eq:beta}) is beneficial for theoretical
purposes but it is computationally inefficient. We now show how the
calculation can be done in a recursive and faster way. The key
point is to find ``effective'' weight vectors $\widetilde \w_i$ such that for every $i$
\begin{eqnarray}
\label{eq:tildew}
\t_i&=& \X_i \w_i =\X \widetilde \w_i\,.
\end{eqnarray}
This can be done by exploiting the fact that $\R$
is bidiagonal.
\begin{proposition}
\label{pro:tildew} The effective weight vectors $\widetilde \w_i$
defined in (\ref{eq:tildew}) and the regression vectors of penalized PLS are determined  by setting $\widetilde
\w_0=\bf{0}$ and $\widehat \bbeta^{(0)}=\bf{0}$ and computing
iteratively
\begin{eqnarray*}
\widetilde \w_i&=&  \w_{i} - \frac{\widetilde \w_{i-1} ^T \X^T \X
\w_{i}}{\widetilde \w_{i-1} ^T \X^T \X \widetilde \w_{i-1}}
\widetilde \w_{i-1}\,,\\
\widehat \bbeta ^{(i)}&=& \widehat \bbeta^{(i-1)} +  \frac{
\widetilde \w_{i}^T \X^T \Y}{\widetilde \w_{i} ^T \X^T \X \widetilde
\w_{i}} \widetilde \w_{i}\,.\\
\end{eqnarray*}
\end{proposition}
 The proof can be found in
the appendix. Combining this result with the PLS algorithm
\ref{algo:pls1}, we obtain the penalized PLS algorithm
\ref{algo:penpls}.

\begin{algorithm}
\caption{Penalized PLS}
\begin{algorithmic}
\label{algo:penpls} \STATE{$\X_1=\X$, $\y$, number of components $m$,
penalty matrix $\P$ (input)} \STATE{ $\M=\left(\I_p
+\P\right)^{-1}$, $\widetilde \w_0={\bf 0}$, $\widehat
\bbeta^{(0)}={\bf 0}$ (initialization)}
 \FOR{i=1,\ldots,m}
\STATE{$\w_i= \M \X_i ^T \y$ (weight vector)}

\STATE{$\widetilde \w_i=  \w_{i} -   \frac{\widetilde \w_{i-1} ^t
\X^T \X \w_{i}}{\widetilde \w_{i-1} ^t \X^T \X \widetilde \w_{i-1}}
\widetilde \w_{i-1}$ (effective weight vector)} \STATE{$\widehat
\bbeta^{(i)}=\widehat \bbeta^{(i-1)} +  \frac{ \widetilde \w_{i}^T
\X^T \Y}{\widetilde \w_{i} ^T \X^T \X \widetilde \w_{i}} \widetilde
\w_{i}$ (regression vector)}

 \STATE{ $\t_i= \X_i
\w_i$ (component)} \STATE{ $\X_{i+1}=\X_i - \mathcal{P}_{\t_i} \X_i$
(deflation)} \ENDFOR
\end{algorithmic}
\end{algorithm}

\subsection{Partial Least Squares and Krylov Subspaces}
It is well-known that PLS is closely connected to Krylov subspaces and
conjugate gradient methods. Quite generally, linear regression
problems can be transformed into algebraic problems in the following
way. The OLS estimator is the solution of the minimization problem
\begin{eqnarray}
\label{eq:ls} \min_{\bbeta} && \|\y - \X \bbeta\|^2\,.
\end{eqnarray}
This is equivalent to finding the solution of the associated normal equation
\begin{eqnarray}
\label{eq:normal}
\A \bbeta&=& \b
\end{eqnarray}
with $\b= \X^T \y $ and $\A= \X^T \X \,$. If the matrix $\A$ is
invertible, the solution of the normal equations is the OLS
estimator $\widehat \bbeta=\A^{-1}\b$. If $\A$ is singular, the solution of
(\ref{eq:normal}) with minimal Euclidean norm is $\A^-\b$. We already mentioned in
Section \ref{sec:pls} that  in the case of high dimensional data,
the matrix $\A$ is often (almost) singular and that the OLS
estimator performs poorly on new data sets. A popular strategy is to
regularize the least squares criterion (\ref{eq:ls}) in the hope of
improving the performance of the estimator. This corresponds to
finding approximate solutions of (\ref{eq:normal}). For example,
Ridge Regression corresponds to the solution of the modified normal
equations
\begin{eqnarray*} \left(\A + \lambda \I_p\right)\bbeta&=& \b\,.
\end{eqnarray*}
Here $\lambda>0$ is the Ridge parameter.  Principal Components
Regression uses the eigen decomposition of $\A$
\begin{eqnarray*}
\A&=& \U \bLambda \U^T=\sum_{i=1} ^p \lambda_i \u_i \u_i ^T
\end{eqnarray*}
and approximates $\A$ and $\b$ via the first $m$ eigenvectors
\[
\begin{array}{rclcrcl}
\A&\approx& \sum_{i=1} ^m \lambda_i \u_i \u_i
^T&,&\b&\approx&\sum_{i=1} ^m \left(\u_i ^T \b\right) \u_i\,.
\end{array}
\]
  It can be shown that the PLS estimators are equal to the approximate solutions of the conjugate gradient method \cite{Hestenes5201}.
This is a procedure that iteratively computes
 approximate solutions of (\ref{eq:normal}) by minimizing the quadratic
 function
\begin{eqnarray}
\label{eq:cg} \phi(\bbeta)=\frac{1}{2} {\bbeta}^T {\A} \bbeta - {
\bbeta}^T {\b}= \frac{1}{2} \left \langle \bbeta,\A \bbeta\right
\rangle -\left \langle \bbeta,\b \right \rangle
\end{eqnarray}
along directions that are ${\bf A}$-orthogonal. The approximate
solution obtained after $m$ steps is equal to the PLS estimator
obtained after $m$ iterations.

The conjugate gradient algorithm is in turn closely related to Krylov
subspaces and the Lanczos
algorithm \cite{Lanczos5001}. The latter is  a method for approximating
eigenvalues. The connection between PLS and these methods is well-elaborated in
\citeasnoun{Phatak0301}. We now
establish a similar connection between penalized PLS and the above mentioned
methods.
Set
\[
\begin{array}{rclcrcl}
\A_{\M}&=& \M \A &\text{ and }& \b_{\M}&=& \M \b\,.
\end{array}
\]
Recall that $\M$ is a symmetric and positive definite matrix that is
determined by the penalty term $\P$. We now illustrate that penalized PLS finds approximate solutions of the
preconditioned normal equation
\begin{eqnarray}
\label{eq:precg}
\A_{\M} \bbeta&=& \b_{\M}\,.
\end{eqnarray}
\begin{lemma}
\label{lemmakrylov}
The space spanned by the weight vectors $\w_1,\ldots,\w_m$  of
penalized PLS  is the same as the space
spanned by the Krylov sequence
\begin{eqnarray}
\label{eq:krylov}
\b_{\M},\A_{\M}  \b_{\M},\ldots,\A_{\M} ^{m-1}  \b_{\M}\,.
\end{eqnarray}
\end{lemma}
This is the generalization of a result for ordinary PLS and can be proven
via induction. Details are given in the appendix.  We denote by
\begin{eqnarray*}
\mathcal{K}^{(m)} &=& \mathcal{K}^{(m)}\left(\A_{\M},\b_{\M}\right)
\end{eqnarray*}
the space that is spanned by the Krylov sequence (\ref{eq:krylov}). This space is
called a Krylov space.
\begin{corollary}
\label{corkrylov}
The penalized PLS  estimator is the solution of the optimization problem
\begin{eqnarray*}
\min& &\|\y-\X\bbeta\|^2\\
\text{subject to}& &\bbeta \in \mathcal{K}^{(m)}.
\end{eqnarray*}
\end{corollary}
\begin{proof}
This follows immediately from proposition \ref{regvector} and the fact that the weight vectors span the Krylov space $\mathcal{K}^{(m)}$.
\end{proof}
We now present the conjugate gradient method for the equation
\begin{eqnarray}
\label{eq:presystem} \A_{\M} \bbeta &=& \b_{\M}\,.
\end{eqnarray}
The Conjugate gradient method is normally applied if the involved
matrix is symmetric. Note that in general, the matrix $\A_{\M}$ is
not symmetric with respect to the canonical inner product, but with
respect to the inner product
\begin{eqnarray*}
\langle \x,\widetilde \x\rangle_{\M^{-1}}&=& \x^T \M^{-1} \widetilde
\x
\end{eqnarray*}
defined by $\M^{-1}$. We can rewrite the quadratic function $\phi$
defined in (\ref{eq:cg}) as
\begin{displaymath}
\phi(\bbeta)= \frac{1}{2} \left \langle \bbeta,\A_{\M} \bbeta\right
\rangle_{\M^{-1}} -\left \langle \bbeta,\b_{\M} \right
\rangle_{\M^{-1}}\,.
\end{displaymath}
We replace the canonical inner product by the inner product defined
by $\M^{-1}$ and minimize this function iteratively along directions
that are $\A_{\M}$-orthogonal.

We start with an initial guess $\bbeta_0={\bf 0}$ and define
$\d_0=\r_0=\b_{\M} - \A_{\M} \bbeta_0 = \b_{\M}$. The quantity
$\d_m$ is the search direction  and $\r_m$ is the residual. For a
given direction $\d_m$, we have to determine  the optimal step size,
that is we have to find
\begin{displaymath}
a_m=\text{arg}\min_{a} \phi\left( \bbeta_m +a \d_m \right) \,.
\end{displaymath}
It is straightforward to check that
\begin{displaymath}
a_m= \frac{\left \langle \d_m,\r_m\right \rangle_{\M^{-1}}}{ \left
\langle \d_m,\A_{\M} \d_{m}\right \rangle_{\M^{-1}}}\,.
\end{displaymath}
The new approximate solution is then
\begin{displaymath}
\bbeta_{m+1}= \bbeta_m + a_m \d_m\,.
\end{displaymath}
After updating the residuals via
\begin{displaymath}
\r_{m+1}= \b_{\M} - \A_{\M} \bbeta_{m+1},
\end{displaymath}
we define a new search direction $\d_{m+1}$ that is
$\A_{M}$-orthogonal to the previous search directions. This is
ensured by projecting the residual $\r_m$ onto the space that is
$\A_{\M}$-orthogonal  to $\d_0,\ldots,\d_m$. We obtain
\begin{displaymath}
\d_{m+1}= \r_{m+1} - \sum_{i=0} ^m \frac{\left \langle
\r_{m+1},\A_{\M} \d_i \right \rangle_{\M^-1}}{\left \langle
\d_{i},\A_{\M} \d_i \right \rangle_{\M^-1} } \d_{i}\,.
\end{displaymath}
\begin{thm}
\label{thm:cg} The penalized PLS algorithm is equal to the conjugate
gradient algorithm for the preconditioned system (\ref{eq:precg}).
\end{thm}
The presentation of the conjugate gradient method above and the
proof of its equivalence to penalized PLS are an extension of the
corresponding results for PLS that is given in \citeasnoun{Phatak0301}.
The proof can be found in the appendix.
Note that there is a different notion of conjugate gradients for
preconditioned systems \cite{Golub8301}. We
transform the preconditioned equation (\ref{eq:cg}) by postmultiplying with
$\M$:
\[
\begin{array}{rclcrcl}
\M \A \M \tilde \bbeta&=& \M\b &\text{ with }& \tilde \bbeta&=& \M^{-1} \bbeta\,.
\end{array}
\]
As the matrix $\M \A \M$ is symmetric, we can apply the ordinary conjugate
gradient algorithm to this equation. This approach differs from the one
described above.
\begin{proposition}
Suppose that $\A=\X^T \X$ is regular. After at most $p$ iterations, the penalized PLS estimator
equals the OLS estimator.
\end{proposition}
\begin{proof}
Using (\ref{eq:restr.ols}), the above statement  is equivalent
to showing that
\begin{eqnarray*}
\widehat \bbeta_{OLS} &\in& \mathcal{K}^{(p)}\,.
\end{eqnarray*}
Hence, we have to show that there is a polynomial $\pi$ of degree $\leq p-1$
such that $\widehat \bbeta_{OLS}=\pi \left(\A_{\M}\right)  \b_{\M}$. As $\M$ is invertible, the OLS estimator is
\begin{displaymath}
\widehat \bbeta_{OLS}= \A^{-1} \b = \A^{-1} \M^{-1} \cdot \M \b = \left(\M
\A\right)^{-1} \M \b =\A_{\M} ^{-1} \b_{\M}\,.
\end{displaymath}
As $\A_M$ is the product of two symmetric matrices and $\M$ is positive
definite, $\A_M$ has a real eigendecomposition,
\begin{eqnarray*}
\A_{\M}&=&\U \Gamma \U^{-1}\,.
\end{eqnarray*}
We define the polynomial $\pi$ via the at most $p$ equations
\begin{eqnarray*}
\pi(\gamma_i)&=&\frac{1}{\gamma_i}\,.
\end{eqnarray*}
It follows immediately that $\pi(\A_{\M})=\A_{\M} ^{-1}$. This concludes the proof.
\end{proof}
\subsection{Kernel Penalized Partial Least Squares}
\label{subsec:kernel}
The computation of the penalized PLS estimator as presented in
algorithm \ref{algo:penpls} involves matrices and vectors of
dimension $p \times p$ and $p$ respectively. If the number of
predictors $p$ is very large, this leads to high computational
costs. In this subsection, we show that we can represent this algorithm in terms of matrices and vectors of dimension $n \times n$
and $n$ respectively.
Let us define the $n\times n$ matrix $\K_{\M}$ via
\begin{displaymath}
\K_{\M}=\left(\langle \x_i, \x_j
\rangle_{\M}\right)=\X \M \X^T\,.
\end{displaymath}
This matrix is called the Gram matrix or the kernel matrix of $\X$.
We conclude from corollary \ref{corkrylov} that the penalized PLS estimator
obtained after $m$ steps is an element of the Krylov space
$\mathcal{K}^{(m)}(\A_{\M},\b_{\M})$. It follows that we can represent
the penalized PLS estimator as
\[
\begin{array}{rclcrcl}
\widehat \bbeta^{(m)}&=& \M \X^T \balpha^{(m)}&,& \balpha^{(m)} &\in &\mathcal{K}^{(m)} \left(\K_{\M},\y\right)\,.
\end{array}
\]
Here, the Krylov space $\mathcal{K}^{(m)} \left(\K_{\M},\y\right)$ is the space
spanned by the vectors
\begin{displaymath}
\y,\K_{\M}\y,\ldots, \K_{\M} ^{m-1} \y\,.
\end{displaymath}
Analogously, we can represent the effective weight vectors by
\[
\begin{array}{rclcrcl}
\widetilde \w_m&=& \M \X^T \widetilde \balpha_m&,& \widetilde \balpha_m &\in &\mathcal{K}^{(m)} \left(\K_{\M},\y\right)\,.
\end{array}
\]
It follows from the definition of the deflation step that
\begin{displaymath}
\X_m ^T \y= \X^T \left(\I_n - \PP_{\t_1,\ldots,\t_{m-1}}\right) \y=
\X^t \left(\y - \widehat \y^{(m-1)}\right)\,.
\end{displaymath}
We conclude that the weight vector $\w_i$ is simply
\[
\begin{array}{rclcrcl}
\w_m&=& \M \X^T \y_{res} ^{(m)}&,&\y_{res}^{(m)}&=&\y - \widehat
\y^{(m-1)}\,.
\end{array}
\]
If we plug in these representations into the penalized PLS algorithm 
\ref{algo:penpls}, we obtain algorithm \ref{algo:kernelpls} that
depends only on $\K_{\M}$ and $\y$.
\begin{algorithm}
\caption{Kernel penalized PLS}
\begin{algorithmic}
\label{algo:kernelpls} \STATE{$\X$, $\y$, number of
components $m$, penalty term $\P$ (input)} \STATE{$\M=\left(\I_p +
\P\right)^{-1}$, $\K_{\M}= \X \M \X^T$, $\balpha^{(0)}=\widetilde
\balpha^{(m)}={\bf 0}$ (initialization)}
 \FOR{i=1,\ldots,m}

  \STATE{$\y_{res}^{(m)}=\y - \widehat \y^{(m-1)}$ (residuals)
   \STATE{$\widetilde \balpha_m=  \y^{(m)} _{res} -   \frac{ \widetilde \balpha_{m-1}
^T \K_{\M}^2 \y^{(m)} _{res}}{\widetilde \balpha_{m-1} ^T  \K_{\M}
^2  \widetilde \balpha_{m-1}}  \widetilde \balpha_{m-1}$ (effective
weight vector)}

\STATE{$\balpha^{(m)}= \balpha^{(m-1)} +  \frac{  \widetilde
\balpha_{m} ^T \K_{\M} \y}{ \widetilde \balpha_{m} ^T \K_{\M} ^2
\widetilde \balpha_{m} } \widetilde \balpha_{m}$ (regression
vector)} \STATE{$\t_i = \K_{\M} \widetilde \balpha_{m}$ (component)}
 \STATE{$\widehat \y ^{(m+1)} = \widehat \y ^{(m)} +  \mathcal{P}_{\t_i} \y$ (estimation of
 $\y$)}}

\ENDFOR
\end{algorithmic}
\end{algorithm}

A kernel version of PLS has already been defined in
\citeasnoun{Raennar9401} in order to speed up the computation of PLS. We
repeat that the speed of the kernel version of penalized PLS does
not depend on the number of predictor variables at all but on the
number of observations. This implies that -- from an algorithmic point
of view -- there are no restrictions in terms of the number of
predictor variables. The importance of this so-called ``dual''
representation also becomes apparent if we want to extend PLS to
nonlinear problems by using the kernel trick. In this paper, the
kernel trick appears in two different versions.

Let us only consider the case of ordinary PLS on B-Splines
transformed variables. Recall that in (\ref{eq:Z}), we transform the
original data $\X$ using a nonlinear function $\Phi$ defined by the
B-Splines.  As algorithm \ref{algo:kernelpls} only relies on inner
products between observations, the nonlinear transformation does not
increase the computational costs. We only have to compute the kernel
matrix of inner  products
\begin{eqnarray*}
\K&=& \left( \left\langle \Phi(\x_i),\Phi(\x_j) \right\rangle  \right)_{i,j=1,\ldots,n}\,.
\end{eqnarray*}
This implies that we do not have to map the data points explicitly using a
function $\Phi$. It suffices to compute the function
\begin{eqnarray}
\label{kernel}
k(\x,\widetilde \x) & = &\left\langle \Phi (\x),\Phi(\widetilde \x) \right \rangle\,.
\end{eqnarray}
The function $k$ is called a kernel. The replacement of the usual inner product by kernel is known as the kernel trick and has turned up to be very
popular in the machine learning community. Instead of defining a nonlinear
map $\Phi$, we define a ``valid'' kernel function $k(x,z)$. E.g., polynomial
relationships can be modeled via kernels of the form
\begin{eqnarray*}
k_d(x,z)&=& \left( 1+ \langle x,z\rangle\right)^d\,,d\in \mathbb{N}\,.
\end{eqnarray*}
Furthermore, it is possible to define kernels for complex data
structures as graphs or text. Literature on the kernel trick and its
applications is abundant. A detailed treatise of  the subject can be
found in \citeasnoun{Schoelkopf0201}. A nonlinear version of PLS using the
kernel trick is presented in \citeasnoun{Rosipal0101}.

If we represent penalized PLS in terms of the kernel matrix $\K_{\M}$,  we realize that penalized
PLS is closely connected to the kernel trick in other respects. Using algorithm \ref{algo:kernelpls} or the definition of the kernel
matrix $\K_{\M}$, we realize that penalized PLS equals  ordinary PLS
with the canonical inner product replaced by the inner product
\begin{displaymath}
\langle \x,\z\rangle_{\M}= \x^T \M \z\,.
\end{displaymath}
This function is called a linear kernel. Why is this a sensible
inner product? Let us consider the eigendecomposition of the penalty
matrix, $\P=\S \bTheta \S^T$. We prefer direction $\s$ such that
$\s^T\P \s$ is small, that is we prefer directions that are defined
by eigenvectors $\s_i$ of $\P$ with a small corresponding eigenvalue
$\theta_i$. If we represent the vectors $\x$ and $\z$ in terms of
the eigenvectors of $\P$,
\[
\begin{array}{ccc}
\widetilde \x= \S^T \x&,& \widetilde \z=\S^T \z\,,
\end{array}
\]
we conclude that
\begin{displaymath}
\langle \x,\z\rangle_{\M}= \widetilde \x ^T \left(\I_p +\bTheta\right)^{-1}
\widetilde \z = \sum_{i=1} ^p \frac{1}{1+\theta_i} \widetilde \x_i \widetilde \z_i
\end{displaymath}
This implies that directions $\s_i$ with  a small eigenvalue $\theta_i$ receive a higher
weighting than directions with a large eigenvalue.
\section{Example: Birth Data}
\label{sec:birth}
\label{sec:example}
In this section, we analyze a real data set describing pregnancy and
delivery for $42$ infants who are sent to a neonatal intensive care
unit after birth. The data are taken from the {\tt R} \cite{R} software
package {\tt exactmaxsel} and are introduced in  \citeasnoun{Boulesteix0602}.
Our goal is to predict the number of days spent in the neonatal
intensive care unit (y) based on the following predictors: birth
weight (in g), birth height (in cm),
 head circumference (in cm), term (in week), age of the mother (in year),
 weight of the mother before pregnancy (in kg), weight of the mother before
 delivery (in kg), height of the mother (in cm), time (in month).
 Some of the predictors are expected to be strongly associated with the
 response (e.g., birth weight, term), in contrast to poor predictors like time or
 height of the mother.

The parameter settings are as follows. We make the simplifying
assumption that $\lambda=\lambda_1=\ldots=\lambda_p$, which reduces
the problem of selecting the optimal smoothing parameter to a
one-dimensional problem. As already mentioned above, we use cubic
splines. Furthermore, the order of difference of adjacent weights is
set to $2$.
\begin{figure}[ht]
{\par\centering\resizebox*{14cm}{8cm}{\rotatebox{270}{{\includegraphics{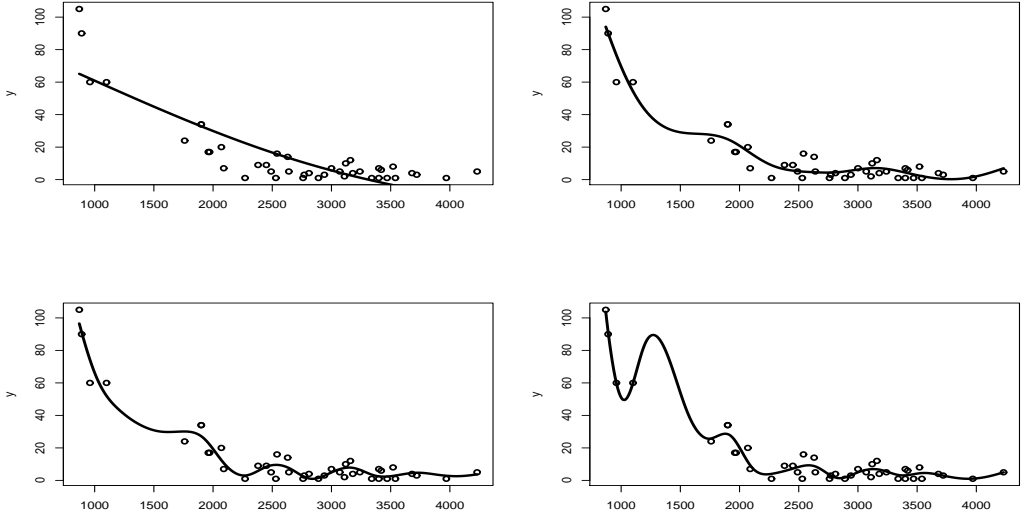}}}}\par}
\caption{Fitted function for the predictor variable ``weight'' using
penalized PLS. The value of $\lambda$ is 2000 and the numbers of
components are $1,5$ (top) and $9,13$ (bottom).} \label{fig:weight}
\end{figure}
The shape of the fitted functions $f_j$ depends on the two model
parameters $\lambda$ and $m$.  We first illustrate that the number
$m$ of penalized PLS components controls the smoothness of the estimated
functions. For this purpose, we  only consider the predictor
variable ``weight''. Figure \ref{fig:weight} displays the fitted
functions obtained by penalized PLS for $\lambda=2000$ and 4
different numbers of  components $m=1,5,9,13$. For small values of $m$, the obtained functions are smooth. For
higher values of $m$, the functions adapt themselves more and more
to the data which leads to overfitting for high values of $m$.

We compare our novel method to PLS without penalization as
described in \cite{Dur2001} and the \texttt{gam()} package in
\texttt{R}. This is the implementation of an adaptive selection
procedure for the basis functions in (\ref{eq:expansion}). More
details can be found in \citeasnoun{Wood0001} and \citeasnoun{Wood0601}. This is the
standard tool for estimating generalized additive models.  The
optimal parameter values of (penalized) PLS are determined by
computing the leave-one-out squared error.  We remark that the split into
training and test set is done before transforming the original
predictors using B-splines. In order to have comparable results, we
normalize the response such that $\text{var}(\y)=1$. The results are summarized in Table
\ref{tab:results2}. Penalized PLS is the best out of the three
method. In particular, it receives a considerably lower error than
PLS without penalization.
\begin{table}[ht]
\vspace{0.2cm} {\centering \begin{tabular}{r|c|c|c}
&leave-one-out-error&$m_{opt}$ & $\lambda_{opt}$\\
\hline PLS &0.159&8 & --\\
penalized PLS & {\bf{0.090}}& 2& 330 \\
GAM& 0.115& -- & --
\end{tabular}
\caption{Optimal model parameters and leave-one-out error for the
birth data set and normalized response.} \label{tab:results2}
\vspace{0.2cm}}
\end{table}
\section{Example:Polymer Data}
\label{sec:polymer}
This data set consists of $p=10$ predictor variables and four
response variables. The number of observations is $n=61$. The data are
taken from a polymer test plant. It  can be downloaded from
\url{ftp://ftp.cis.upenn.edu/pub/ungar/chemdata/}. The predictor
variables are measurements of controlled variables in a polymer
processing plant (e.g. temperatures, feed rates ...). No more
details on the variables are given due to confidentiality reasons.
As in the last section, we first scale each response variable to
have a variance equal to $1$. Again, we compare penalized PLS to
PLS and \texttt{gam()}. The results are summarized in Table
\ref{tab:results3}.
\begin{table}[ht]
\vspace{0.2cm} {\centering \begin{tabular}{r|c|c|c|c} &
$1^{\text{st}}$ response & $2^{\text{nd}}$ response &
$3^{\text{rd}}$  response & $4^{\text{th}}$
response\\
\hline PLS &0.672 & 0.863 & 0.254 & 0.204\\
penalized PLS & 0.607 & {\bf{0.801}} & {\bf{0.206}} & {\bf{0.164}} \\
GAM& {\bf{0.599}} & 0.881 & 0.218 & 0.182
\end{tabular}
\caption{Leave-one-out error for the polymer data set and normalized
response.} \label{tab:results3} \vspace{0.2cm}}
\end{table}
For all four response variables,  penalized PLS is better than PLS
without penalization. Penalized PLS is also better that GAM for
three out of the four response variables, although the difference is
considerably smaller.
\section{Concluding Remarks}
\label{sec:conclusion}
In this work, we proposed an extension of Partial Least Squares
Regression using penalization techniques. Apart from its
computational efficiency (it is virtually as fast as PLS), it also
shares a lot of mathematical properties of PLS. Our novel 
method  obtains good results in applications. In the two
examples that are discussed, penalized PLS clearly outperforms PLS without
penalization. Furthermore, the results indicate that it is a
competitor of \texttt{gam()}  in the case of very high-dimensional data.

We might think of other penalty terms. \citeasnoun{Kondylis0601} consider a
preconditioned version of PLS by giving weights to the predictor
variables. Higher weights are given to those predictor variables that are
highly correlated to the response. These weights can be expressed in terms
of a penalty term. \citeasnoun{Goutis9602} combine PLS with an additive penalty
term to data derived from near infra red spectroscopy. The penalty term
controls the smoothness of the regression vector. 

The introduction of a penalty term can easily be adapted to other dimension
reduction techniques. For example for Principal Components Analysis, the
penalized optimization criterion is
\begin{displaymath}
\max_{\w} \frac{\text{var} (\X \w)}{\w^T \w + \w^T \P \w}\,.
\end{displaymath}
PLS can handle multivariate responses $\Y$. The natural extension of criterion
(\ref{eq:crit1}) is the following.
\begin{displaymath}
\max_{\w} \frac{\left \|\text{cov}( \X \w,\Y)\right \|^2}{\w^T \w}=
\max_{\w} \frac{\w^T \X^T \Y \Y^T \X \w}{\w^T \w}\,.
\end{displaymath}
Using Lagrangian multipliers, we deduce that the solution is the eigenvector
of the matrix
\begin{displaymath}
\B=\X^T \Y \Y^T \X
\end{displaymath}
that corresponds to the largest eigenvalue of $\B$. This eigenvector is
usually computed in an iterative fashion. If we want to apply penalized PLS
for multivariate responses, we compute
\begin{displaymath}
\max_{\w} \frac{\w^T \X^T \Y \Y^T \X \w}{\w^T \w +\w^T \P\w}\,.
\end{displaymath}
The solution fulfills
\begin{displaymath}
\B \w = \gamma \left(\I_p + \P\right) \w, \, \gamma \in \mathbb{R}\,.
\end{displaymath}
This is called a generalized eigenvalue problem or a matrix pencil. Note
that for multivariate $\Y$, the equivalence of SIMPLS and NIPALS does not
hold, so we expect the penalized versions of these methods to be different
as well. There are kernel versions for PLS with multivariate $\Y$
\cite{Raennar9401,Rosipal0101}, hence we can also represent multivariate
penalized PLS in terms of kernel matrices.
\section*{Acknowledgement}
This research was supported by the Deutsche Forschungsgemeinschaft (SFB 386,
``Statistical Analysis of Discrete Structures'').

\bibliographystyle{agsm}

\appendix
\section{Proofs}
We recall that for $k<i$
\begin{eqnarray}
\label{defl1} \X_i= \prod_{j=k} ^{i-1} \left(\I_n
-\mathcal{P}_{\t_j}\right) \X_k = \left(\I_n
-\mathcal{P}_{\t_k,\ldots,\t_{i-1}}\right) \X_k
\end{eqnarray}
The last equality follows from the fact that the components $\t_i$
are mutually orthogonal. In particular, we obtain
\begin{eqnarray}
\label{defl2} \X_i&=& \left(\I_n
-\mathcal{P}_{\t_1,\ldots,\t_{i-1}}\right) \X\,.
\end{eqnarray}

\begin{proof}[Proof of lemma \ref{lemma:R}]
First note that (\ref{defl2}) is equivalent to $\X= \X_j
+\mathcal{P}_{\t_1,\ldots,\t_{j-1}}\X\,$. It follows that
\begin{eqnarray}
\label{rij}
\X \w_j&=& \X_j \w_j + \mathcal{P}_{\t_1,\ldots,\t_{j-1}}\X\ \w_j = \t_j + \sum_{i=1} ^{j-1} \frac{\t_i ^T \X \w_j}{\t_i ^T \t_i} \t_i\,.
\end{eqnarray}
As all components $\t_i$ are mutually orthogonal,
\begin{eqnarray*}
\t_i ^T \X \w_j&=& \begin{cases}
\t_i ^T \t_i \not =0 &,i=j\\
0&, i>j\\
*&,\text{otherwise}
\end{cases}\,.
\end{eqnarray*}
We conclude that $\R$  is an upper triangular matrix with all diagonal
elements $\not=0$. Furthermore, it follows from (\ref{rij}) that all vectors
$\X\w_j$ are linear combinations of the components $\t_1,\ldots,\t_j$. This
implies that the columns of $\X\W$ and the columns of $\T$ span the same
space. Finally, we have to show that $\R$ is bidiagonal. To prove this, we
show that $\X_i \w_j=0$ for  $j<i\,$.  The condition $i>j$ implies
(recall (\ref{defl1})) that $\X_i= \X_j -
\mathcal{P}_{\t_j,\ldots,\t_{i-1}} \X_j$ and consequently
\begin{eqnarray*}
\X_i \w_j&=& \X_j \w_j - \mathcal{P}_{\t_1,\ldots,\t_{i-1}} \X_j \w_j= \t_j -  \mathcal{P}_{\t_1,\ldots,\t_{i-1}} \t_j
\stackrel{j\leq i-1}{=} \t_j - \t_j={\bf{0}}\,.
\end{eqnarray*}
This implies that for $i-1>j$
\begin{eqnarray*}
\t_i ^T \X \w_j&=& \t_i ^T \left(\X_{i} +\mathcal{P}_{\t_1,\ldots,\t_{i-1}} \X
\right) \w_j \\
&=& \t_i ^T \left( \X_i \w_j +\mathcal{P}_{\t_1,\ldots,\t_{i-1}} \X \w_j \right)\\
&=& \t_i ^T \left(\mathcal{P}_{\t_1,\ldots,\t_{i-1}} \X \w_j \right)={\bf{0}}\,.
\end{eqnarray*}
\end{proof}
\begin{proof}[Proof of proposition \ref{pro:tildew}]
For $i=1$, we have $\tilde \w_1= \w_1$ as $\X_1=\X$. For a general $i$, we
have
\begin{displaymath}
\t_{i+1}= \X_{i+1} \w_{i+1} =\left(\X - \mathcal{P}_{\t_1,\ldots,\t_i} \X \right) \w_{i+1}= \X \w_{i+1} - \mathcal{P}_{\t_{i}} \X \w_{i+1}\,.
\end{displaymath}
The last equality holds as $\R=\T^T \X\W$ is bidiagonal. Using  formula (\ref{eq:proj})
for the projection operator, it follows that
\begin{eqnarray*}
\t_{i+1}&=& \X \w_{i+1} - \X \frac{\tilde \w_{i} \tilde \w_i ^T}{\tilde
\w_{i} ^T \X^T \X \tilde \w_{i}} \X^t  \X \w_{i+1}\,.
\end{eqnarray*}
We conclude that
\begin{eqnarray*}
\tilde \w_{i+1}&=& \w_{i+1} -   \frac{\tilde \w_i ^T \X^T \X \w_{i+1}}{\tilde \w_{i} ^T
\X^T \X \tilde \w_{i}} \tilde \w_i\,.
\end{eqnarray*}
The regression estimate after $i+1$ steps is
\begin{eqnarray*}
\X \widehat \bbeta^{(i+1)}&=& \mathcal{P}_{\t_1,\ldots,\t_{i+1}} \Y\\
&=& \X \widehat \bbeta^{(i)} + \mathcal{P}_{\t_{i+1}} \Y\\
&=& \X \widehat \bbeta^{(i)} + \X \frac{ \tilde \w_{i+1} \tilde
\w_{i+1}^T}{\tilde \w_{i+1} ^T \X^T \X \tilde \w_{i+1}} \X^T \Y\,.
\end{eqnarray*}
This concludes the proof.
\end{proof}
\begin{proof}[Proof of lemma \ref{lemmakrylov}]
We use induction. For $m=1$ we know that $\w_1= \b_{\M}$. For a
fixed $m>1$, we conclude from the induction hypothesis and
lemma \ref{lemma:R} that every
vector $\s$ that lies in the span of $\t_1,\ldots,\t_m$ is of the form
\begin{eqnarray}
\label{eq:s}
\s= \X {\v} &,& \v \in \text{span}\{\w_1,\ldots,\w_m\}=\mathcal{K}^{(m)}\,.
\end{eqnarray}
We conclude that
\begin{displaymath}
\X_{m+1} ^T \y= \left(\X -  \mathcal{P}_{\t_1,\ldots,\t_m} \X \right)^T
\y= \X^T \y - \X^T \mathcal{P}_{\t_1,\ldots,\t_m} \y \stackrel{(\ref{eq:s})}{=} \b - \X^T \X \s
\end{displaymath}
and that
\begin{displaymath}
\w_{m+1}= \M \X_{m+1} ^T \y= \M \b - \M \A \s =\b_{\M} - \A_{\M} \s \in \mathcal{K}^{(m+1)}\,.
\end{displaymath}
\end{proof}
In the rest of the appendix, we show the equivalence of penalized
PLS and the preconditioned conjugate gradient method.
\begin{lemma}
\label{lem:krylov} We have
\begin{displaymath}
\text{span}\left \{ \d_0,\ldots,\d_{m-1}\right \}= \text{span}\left
\{ \r_0,\ldots,\r_{m-1}\right \} = \text{span}\left \{
\x_1,\ldots,\x_{m}\right \}=\mathcal{K}^{(m)}\,.
\end{displaymath}
\end{lemma}
This can be proven via induction.
\begin{lemma}
\label{lem:xm} We have
\begin{displaymath}
\bbeta_m= \sum_{i=0} ^{m-1} \frac{\left \langle \d_i,\b_{\M}\right
\rangle_{\M^{-1}}}{ \left \langle \d_i,\A_{\M} \d_{i}\right
\rangle_{\M^{-1}}} \d_i
\end{displaymath}
\end{lemma}
\begin{proof}
This corresponds to the iterative definition of $\bbeta_{m+1}$. We
only have to show that
\begin{eqnarray*}
\left \langle \d_i,\r_i\right \rangle_{\M^{-1}} &=& \left \langle
\d_i,\b_{\M}\right \rangle_{\M^{-1}}\,.
\end{eqnarray*}
Note that
\begin{displaymath}
\r_i= \b - \sum_{j=0} ^{i-1} a_j \A_{\M} \d_j\,.
\end{displaymath}
As $\d_i$ is $\A_{\M}$-orthogonal onto all directions
$\d_j\,,\,j<i$, the proof is complete.
\end{proof}
Now we are able to proof the equivalence of penalized PLS and the
conjugate gradient method.
\begin{proof}[Proof of theorem \ref{thm:cg}]
As the search directions $\d_i$ span the Krylov space
$\mathcal{K}^{(m)}$, we can replace the matrix $\W$ in
(\ref{eq:beta}) by the matrix
$\D=\left(\d_0,\ldots,\d_{m-1}\right)$. As the search directions are
$\A_{\M}$-orthogonal, we have
\begin{eqnarray*}
\widehat \bbeta_{PPLS} &=& \D\left(\D^T \A \D\right)^{-1} \D^T \b\\
&=& \D\left(\D^T \M^{-1} \A_{\M} \D\right)^{-1} \D^T \M^{-1}\b_{\M}\\
&=&  \sum_{i=0} ^{m-1} \frac{\left \langle \d_i,\b_{\M}\right
\rangle_{\M^{-1}}}{ \left \langle \d_i,\A_{\M} \d_{i}\right
\rangle_{\M^{-1}}} \d_i
\end{eqnarray*}
and this equals the formula in lemma \ref{lem:xm}.
\end{proof}
\end{document}